\documentclass[10pt,leqno,twoside]{article}
\usepackage{amssymb}
\usepackage{amscd}
\usepackage{amsthm}
\theoremstyle{plain}

\newtheorem{Thm}{Theorem}[section]
\newtheorem{Cor}[Thm]{Corollary}
\newtheorem{Lemma}[Thm]{Lemma}

\newtheorem{Obs}[Thm]{Observation}
\newtheorem{Prop}[Thm]{Proposition}

\theoremstyle{definition}
\newtheorem{Def}[Thm]{Definition}

\newtheorem{Rmk}[Thm]{Remark}

\newcommand{\B}{{\cal B}}
\newcommand{\Cs}{C^{\ast}}

\newcommand{\A}{{\cal A}}
\newcommand{\M}{{\cal M}}

\newcommand{\K}{{\cal K}}

\newcommand{\HH}{{\cal H}}

\newcommand{\R}{{\mathbb R}}
\newcommand{\N}{{\mathbb N}}
\newcommand{\Z}{{\mathbb Z}}
\newcommand{\C}{{\mathbb C}}

\newcommand{\arrow}{\rightarrow}

\newcommand{\cxk}{C_0(X) \otimes {\cal K}}

\newcommand{\Cn}{{\cal C}_n}
\newcommand{\Un}{{\cal U}_n}
\newcommand{\UH}{{\cal U(H)}}
\newcommand{\cxendo}{$C_0(X)$-endomorphism}
\newcommand{\GE}{\Gamma_0(E)}
\newcommand{\otx}{\otimes_X}
\newcommand{\Intn}{{\cal I}_n}
\newcommand{\Va}{V_{\alpha}}

\newcommand{\EnK}{End(\K)}
\newcommand{\EnnK}{End^n(\K)}

\long\def\MSC#1\EndMSC{\def\arg{#1}\ifx\arg\empty\relax\else
     {\par\narrower\noindent%
     2000 Mathematics Subject Classification: #1\par}\fi}

\newcommand{\Proof}{\textbf{Proof:}~}

\begin{document}                                               
       
\title{Endomorphisms of stable continuous-trace $C^*$-algebras}
\author{Ilan Hirshberg}

\date{}
\maketitle

\begin{abstract}
We classify spectrum-preserving endomorphisms of stable continuous-trace 
$\Cs$-algebras up to inner automorphism by a surjective multiplicative 
invariant taking values in finite dimensional vector bundles over the 
spectrum. Specializing to automorphisms, this gives a different approach to 
results of Phillips and Raeburn. 
\end{abstract}

\MSC
46L05, 46M20
\EndMSC

\section{Introduction}

The structure of automorphisms of continuous-trace $\Cs$-algebras was 
studied by Lance (\cite{Lance}) and Smith (\cite{Smith}) (the case of 
$C(X)-$automorphisms of $C(X)\otimes \K$, or more precisely 
$C(X)\otimes \B(\HH)$), and subsequently by Phillips and Raeburn (\cite{PR}, 
see also \cite{Rosenberg},\cite{RW}) for general continuous trace 
$\Cs$-algebras.

This paper is concerned with the question of generalizing some of those 
results to the situation of `unital' endomorphisms of stable 
continuous-trace $\Cs$-algebras. We refer the reader to \cite{RW} and 
\cite{Dixmier} for general references on the theory of continuous-trace 
$\Cs$-algebras.

Throughout the paper, $X$ will denote a locally compact Hausdorff space; 
$\HH$ will denote a separable Hilbert space and $\B(\HH)$ will denote the 
bounded operators on $\HH$; $\M(\A)$ will denote the multiplier algebra of 
$\A$; $\K$ will denote the algebra of compact operators on a separable 
Hilbert space; automorphisms, endomorphisms and homomorphisms between 
$\Cs$-algebras and $\Cs$-bundles will always be $*$-homomorphisms; 
$\delta(\A)$ will denote the Dixmier-Douady invariant of a continuous-trace 
$\Cs$-algebra $\A$. By a unital homomorphism between two non-unital 
$\Cs$-algebras we mean a homomorphism which maps approximate units to 
approximate units.

We recall that every stable continuous-trace $\Cs$-algebra for which all the 
irreducible representations act on a separable space is isomorphic to the 
algebra of sections of a locally trivial $\K$-bundle over a locally compact 
Hausdorff space. Let $E$ be a locally trivial $\K$-bundle over $X$, and let 
$\A = \GE$. Denote by $Aut_X(\A)$ the automorphisms of $\A$ which commute 
with the multiplication action of $C_0(X)$ on $\A$, and by $Inn(\A)$ the 
group of inner automorphisms of $\A$, i.e. those automorphisms of $\A$ which 
are given by $\alpha(A) = UAU^*$ for some $U\in\M(\A)$. Denote by 
$Homeo_{\delta}(X)$ the group of homeomorphisms of $X$ which preserve 
$\delta(A)\in H^3(X;\Z)$, where cohomology here is taken to mean \v{C}ech 
cohomology. The result of Phillips and Raeburn, specialized to this case, 
can be summarized in the following  exact sequence.
\begin{equation}
1 \arrow Inn(\A) \arrow Aut_X(\A) \arrow H^2(X;\Z) \arrow 1
\label{eqn:PRseq1}
\end{equation}
Of course, some of this structure must be lost when generalizing to 
endomorphisms. The fact that the endomorphisms do not form a group means 
that speaking of exact sequences won't make much sense.

This paper will be devoted to providing an adequate generalization 
to sequence (\ref{eqn:PRseq1}). Our invariant takes values in vector bundles over $X$. We will 
show that the associated vector bundle of a $C_0(X)$-endomorphism 
$\alpha$, denoted $vect(\alpha)$, classifies $\alpha$ up to inner 
automorphism (i.e. $vect(\alpha) \cong vect(\beta) \Longleftrightarrow 
\alpha = AdU \circ \beta$); that all numerable vector bundles over $X$ 
arise in this way; and that we have a product formula $vect(\alpha \circ 
\beta) \cong vect(\alpha) \otimes vect(\beta)$ (in other words, $vect$ is a 
semigroup homomorphism from the semigroup of endomorphisms under composition 
to vector bundles under tensor products).

The invariant for $C_0(X)$-automorphisms in this approach will be 
a complex line bundle over $X$. So if we specialize the above results to 
automorphisms, then we get a similar sequence to (\ref{eqn:PRseq1}) (with the 
group of complex line bundles over $X$, which is isomorphic to $H^2(X,\Z)$, on the right).

The structure of this paper is as follows. In section 2, we study the space 
of endomorphisms of $\K$. In section 3, we discuss the case of \cxendo s of $\cxk$. 
In section 4, we generalize the main 
results of section 3 to algebras with possibly non-trivial Dixmier-Douady 
invariant.

\section{Endomorphisms of $\K$}

Let $\alpha$ be an endomorphism of $\K$ which is unital, in the sense that 
it maps approximate units to approximate units.

Let $n$ be the trace of $\alpha(P)$ for a minimal projection $P\in\K$, where 
by trace we mean the unique (unbounded) trace on $\K_+$ whose value on 
minimal projections is 1. This $n$ clearly is independent of the choice of 
$P$, and will be called $index(\alpha)$.

The following simple lemma and its short proof was conveyed to the author by 
W. B. Arveson.

\begin{Lemma}
There are isometries $S_1,...,S_n \in \M(\K) \cong \B(\HH)$ with
$$
\sum_{i=1}^nS_iS_i^* = 1
$$
such that for all $A\in\K$,
$$
\alpha(A) = \sum_{i=1}^nS_iAS_i^*
$$
\end{Lemma}
\Proof Letting $\K$ act irreducibly on some space $\HH$, we observe 
that $\alpha$ gives a representation of $\K$ on $\HH$. This representation 
decomposes into a direct sum of $n$ irreducible representations on pairwise 
orthogonal subspaces $\HH_1 \oplus ... \oplus \HH_n = \HH$, and there are 
isometries $S_i$ from $\HH$ onto $\HH_i$, $i=1,...,n$ which implement 
$\alpha|_{\HH_i}$, i.e. $S_iAS_i^*\xi_i = \alpha(A)\xi_i$ for 
$\xi_i\in\HH_i$. The lemma follows from this immediately. \qed
\paragraph{}
Notice that we have $S_iA = \alpha(A)S_i$, $i=1,...,n$. Let
$$
\Va = \{T\in\B(\HH) \; | \; TA=\alpha(A)T \; \forall A\in\K \}
$$
$\Va$ is clearly a vector space. It also has a natural inner product on it. 
If $T,S\in\Va$ then $T^*SA=AT^*S$ for all $A\in\K$, and therefore $T^*S$ is 
a scalar multiple of the identity, and indeed $\left < S , T \right > = 
T^*S$ is an inner product.

\begin{Lemma}
The $S_i$'s are an orthonormal basis for $\Va$.
\end{Lemma}
\Proof  Let $T\in\Va$  then $T = (\sum_{i=1}^nS_iS_i^*)T = \sum_{i=1}^n 
\left <S_i,T \right >S_i$, as required. Orthonormality is immediate. \qed
\paragraph{}
If $T_1,...,T_m$ are isometries such that $\alpha(A) = 
\sum_{i=1}^mT_iAT_i^*$ for all $A\in\K$ then clearly $m=index(\alpha)=n$, 
and the $T_i$'s form another orthonormal basis for $\Va$. Conversely, if 
$T_1,...,T_n$ are an orthonormal basis for $\Va$ then it is straightforward 
to verify that $\alpha(A) = \sum_{i=1}^nT_iAT_i^*$ for all $A\in\K$.

Denote by $End(\K)$ the space of unital endomorphisms of $\K$, endowed with 
the point-norm topology, and by $\EnnK$ the subset of endomorphisms of index 
$n$. Notice that $\EnK$ (and hence all the $\EnnK$'s) is metrizable, hence 
paracompact (if $A_1,A_2,...$ is a dense sequence in the unit sphere of 
$\K$, then $d(\alpha,\beta) = \sum_k\|\alpha(A_k)-\beta(A_k)\|/2^k$ is a 
metric inducing the point-norm topology).

\begin{Prop}  \label{prop:continuous n}
$index:End(\K) \arrow \N$ is continuous (so $\EnnK$ is closed and open in 
$End(\K)$ for all $n$).
\end{Prop}
\Proof Let $\alpha_k \arrow \alpha$. Let $P$ be a minimal 
projection, then for sufficiently large $k$, we have 
$\|\alpha_k(P)-\alpha(P)\|<1$. Therefore, $\alpha_k(P)$ and $\alpha(P)$ are 
Murray--von-Neumann equivalent, and in particular they have the same trace, 
i.e. $index(\alpha_k) \arrow index(\alpha)$.
(cf \cite{Price} Proposition 2.3)
\qed
\paragraph{}
We denote by $\Cn \subseteq \B(\HH)^n$ the space of all $n$-tuples of 
isometries satisfying the Cuntz relations, equipped with the strong operator 
topology.  We remark that $\Cn$ is metrizable (since the unit sphere of 
$\B(\HH)$ is metrizable in the SOT), hence paracompact.
We described above a surjective map $\pi:\Cn \arrow \EnnK$, given by 
$\pi(S_1,...,S_n)(A) =  \sum_{i=1}^nS_ixS_i^*$. Suppose 
$(S_1^{(k)},...,S_n^{(k)}) \arrow (S_1,...,S_n)$ in $\Cn$, and let $A=\xi 
\otimes \bar{\eta} \in \K$ a rank 1 operator, then
$$
\pi(S_1^{(k)},...,S_n^{(k)})(A) =  \sum_{i=1}^nS_i^{(k)}AS_i^{(k)*} =
$$ $$ =
\sum_{i=1}^nS_i^{(k)}\xi \otimes \overline{S_i^{(k)*}\eta} \arrow 
\sum_{i=1}^nS_i^*\xi \otimes \overline{S_i^*\eta} = \pi(S_1,...,S_n)(A)
$$ so $\pi$ is continuous.

We define the \emph{universal intertwining bundle}
$$
\begin{CD}
\C^n \\
@VVV \\
\Intn \\
@VVV \\
\EnnK
\end{CD}
$$
to be the sub-bundle of $\EnnK \times \B(\HH)$ whose fiber over $\alpha$ is 
$\Va$ (where the topology on $\B(\HH)$ is the strong operator topology). $$
\begin{CD}
\Un \\
@VVV \\
\Cn \\
@V{\pi}VV \\
\EnnK
\end{CD}
$$
is the corresponding principal $\Un$-bundle (the bundle of orthonormal 
frames).

\begin{Lemma} \label{lemma:local cross sections}
$\pi:\Cn \arrow \EnnK$ has local cross sections. In other words, the bundles 
above are locally trivial.
\end{Lemma}
\Proof (cf \cite{RW}, Proposition 1.6) We think of $\K$ as acting 
irreducibly on $\HH$. Fix $\alpha\in \EnnK$, and let 
$P=\nu\otimes\bar{\nu}\in\K$ be a fixed minimal projection. Let 
$\omega_1,...,\omega_n$ be an orthonormal basis for $\alpha(P)\HH$, and 
define $S_i\in\B(\HH)$, $i=1,..,n$ by $S_i\xi = \alpha(\xi \otimes 
\bar{\nu})\omega_i $.

For any $i,j\in{1,...,n}$, and $\xi,\eta\in\HH$, we have
$$
\left < S_i\xi,S_j\eta \right >
=
\left < \alpha(\xi \otimes \bar{\nu})\omega_i , \alpha(\eta \otimes 
\bar{\nu})\omega_j \right >
=
\left < \alpha(\eta \otimes \bar{\nu})^*\alpha(\xi \otimes   
\bar{\nu})\omega_i,\omega_j \right >
=
$$ $$ =
\left < \alpha((\eta \otimes \bar{\nu})^*\xi \otimes 
\bar{\nu})\omega_i,\omega_j \right >
=
\left < \xi , \eta \right >   \left <\alpha(\nu \otimes \bar{\nu}) 
\omega_i,\omega_j \right>
=
\left < \xi , \eta \right >   \left < \omega_i  , \omega_j  \right >
=
\left < \xi , \eta \right >  \delta_{ij}
$$
so $S_1,...,S_n$ are isometries with orthogonal ranges.

Let $A = \xi \otimes \bar{\eta} $, and let $\zeta \in \HH$, then
$$
S_i A \zeta
=
(S_i \xi \otimes \bar{\eta}) \zeta
=
\left < \zeta , \eta \right > S_i\xi
=
\left < \zeta , \eta \right > \alpha(\xi \otimes \bar{\nu})\omega_i
=
$$ $$ =
\alpha((\xi \otimes \bar{\eta}) (\zeta \otimes \bar{\nu})) \omega_i
=
\alpha(A)\alpha(\zeta \otimes \bar{\nu}))\omega_i
=
\alpha(A)S_i\zeta
$$
so $S_i$ is an intertwiner.

Therefore, $\pi(S_1,...,S_n) = \alpha$. Now, define $f_1,...,f_n : \EnnK 
\arrow \HH$ by
$$
f_i(\beta) = \beta(P)\omega_i
$$
The $f_i$ are all clearly continuous, and $f_i(\alpha) = \omega_i$. 
Therefore, there is an open neighborhood $N$ of $\alpha$ in which 
$f_1,...,f_n$ are linearly independent. Let \linebreak $g_1,...,g_n : N 
\arrow \HH$ be obtained from $f_1,...,f_n$ via the Gram-Schmidt 
orthonormalization process, then the $g_i$'s are continuous as well. Now, 
let \linebreak $s_1,...,s_n : N \arrow \B(\HH)$ be given by
$$
s_i(\beta) \xi = \beta(\xi \otimes \bar{\nu}) g_i(\beta) \; , \; \xi \in \HH
$$
then the $s_i$'s are clearly continuous (when $\B(\HH)$ is given the strong 
operator topology), and from the above, we have $ 
(s_1(\beta),...,s_n(\beta)) \in \Cn$ for all $\beta$, and 
$\pi(s_1(\beta),...,s_n(\beta)) = \beta$, so we have a local cross section, 
as required.
\qed

\begin{Thm} \label{thm:contractibility}
$\Cn$ is contractible.
\end{Thm}
\Proof Let $v_t : L^2(\R_+) \arrow L^2(\R_+)$ be the unilateral 
shifts, $p_t = 1-v_tv_t^*$, $(T_1',...,T_n')$ some fixed $n$-tuple of 
isometries satisfying the Cuntz relations on a separable Hilbert space 
$\HH'$. Identify $\HH = \HH' \otimes L^2(\R_+)$, and let $V_t = 1 \otimes 
v_t$, $P_t = 1 \otimes p_t$, $T_i = T_i' \otimes 1$, $i=1,..,n$. Define a 
family of maps $\Phi_t : \Cn \arrow \Cn$ by \[\Phi_t (S_1,...,S_n) = 
(...,V_tS_iV_t^* + P_tT_iP_t,...)\] then it is easy to see $\Phi_t$ is 
SOT-continuous in $t$, that $\Phi_t (S_1,...,S_n) \in \Cn$ for all 
$(S_1,...,S_n) \in \Cn$, and that and $SOT-\lim\Phi_t(S_1,...,S_n) = 
(T_1,...,T_n)$ for any $(S_1,...,S_n) \in \Cn$, so $\Phi_t$ gives us a 
contraction as required.
\qed
\paragraph{}
Recall that an open cover $\{U_{\lambda}\}$ of a topological space $Y$ is said to be \emph{numerable} if there is a locally finite partition of unity $\{f_{\mu}\}$ such that the open cover $\{f_{\mu}^{-1}((0,1])\}$ refines $\{U_{\lambda}\}$. Notice that if $Y$ is paracompact then any open cover is numerable. A principal $\Un$-bundle $\zeta$ over $Y$ is said to be \emph{numerable} if there is a numerable open cover $\{U_{\lambda}\}$ of $Y$ such that $\zeta|_{U_{\lambda}}$ is trivial for all $\lambda$. If $Y$ is paracompact then any locally trivial bundle is numerable. The pull back of a numerable bundle is also numerable.
A numerable principal $\Un$-bundle over a
space $B$ is said to be \emph{universal} if
\begin{enumerate}
\item Any numerable principal $\Un$-bundle over some space
$Y$ is isomorphic to a pull-back bundle of this bundle via a map $Y \arrow 
B$.
\item Two maps $Y \arrow B$ give rise to isomorphic pull-back bundles if and 
only if they are homotopic.
\end{enumerate}
A universal Hermitian vector bundle is defined in the analogous way, and a 
Hermitian vector bundle is universal if and only if its associated principal 
$\Un$-bundle is universal.

By a theorem of Dold (\cite{Dold} Theorem 7.5), a numerable principal $\Un$-bundle 
is universal if and only if the total space is 
contractible.

\begin{Cor} \label{cor:universal bundle}
The above bundles are universal.
\end{Cor}

\begin{Rmk} For any fixed $(S_1,...,S_n) \in \Cn$ we can define a map $\UH 
\arrow \Cn$ by $U \mapsto (US_1,...US_n)$. This map can readily be seen to 
be a homeomorphism, so in fact all the $\Cn$'s are homeomorphic.
\end{Rmk}

\begin{Rmk} The focus of this paper is ordinary (=complex) $\Cs$-algebras, 
however it is worth noticing that the analogous results of this section 
along with their proofs hold for endomorphisms of the real $\Cs$-algebra of 
compact operators on a real separable Hilbert space as well (replacing 
complex bundles by real bundles, unitaries by orthogonals and so on).
\end{Rmk}

\section{$C_0(X)$-Endomorphisms of $\cxk$}
We now turn to study endomorphisms of $\cxk$ which commute with the 
multiplier action of $C_0(X)$. We refer to those as \emph{\cxendo s}.

\begin{Lemma} \label{lemma:endo as map}
Let $\sigma:X \arrow End(\K)$ be a continuous map, then there is a \cxendo ~ 
$\alpha$ of $\cxk$ given by $\alpha(f)(x) = \sigma(x)(f(x))$. Furthermore, 
any \cxendo ~ $\alpha$ of $\cxk$ is of this form.
\end{Lemma}
\Proof The proof is straightforward (and virtually identical to 
Lemma 4.28 in \cite{RW}, using endomorphisms instead of automorphisms), and 
the details are left to the reader.

\paragraph{}

We henceforth denote this map by $x \mapsto \alpha_x$.

\begin{Def}
$$
vect(\alpha) = \{(x,T) \; | \; TA = \alpha_x(A)T \; \; \forall A\in\K\} 
\subseteq X \times \B(\HH)
$$
\end{Def}

Notice that $vect(\alpha)$ is the pull-back of the universal intertwining 
bundle(s) via the map $x \arrow \alpha_x$ (where possibly the dimension may 
be different on different components of $X$; by Proposition 
\ref{prop:continuous n}, the dimension will be locally constant). In particular, $vect(\alpha)$ is numerable.

\begin{Cor} For any numerable Hermitian vector bundle $v$ there exists 
a \cxendo~ $\alpha$ such that $vect(\alpha) \cong v$. Furthermore, 
$vect(\alpha)\cong vect(\beta)$ if and only if $\alpha$ is connected to 
$\beta$ via a path of \cxendo s.
\end{Cor}

\begin{Thm} \label{thm:product formula}
Let $\alpha$,$\beta$ be \cxendo s, then
$$
vect(\alpha \circ \beta) \cong vect(\alpha)\otimes vect(\beta)
$$
\end{Thm}
\Proof Define a map
$$
vect(\alpha)\otimes vect(\beta) \arrow vect(\alpha \circ \beta)
$$
by
$$
(x,T \otimes S) \mapsto (x,TS)
$$
Verifying that this map is a well defined bundle isomorphism is immediate. 
\qed
\paragraph{}
The associated vector bundle to a \cxendo ~ classifies the \cxendo ~ up to 
homotopy. A more natural equivalence relation on \cxendo s is classification 
up to inner automorphism, i.e., $\alpha \sim \beta$ if there is a unitary 
$U$ in $\M(\cxk)$ such that $Ad\;U \circ \alpha = \beta$. The following 
proposition shows that the two equivalence relations coincide.

\begin{Obs} \label{obs:hom bundle}
$Hom(vect(\alpha),vect(\beta))$ is isomorphic to the sub-bundle of the $X 
\times \B(\HH)$ whose fiber over $x$ is $span\{TS^* \; | \; T \in 
vect(\beta) \; , \; S \in vect(\alpha)\}$.
\end{Obs}

\begin{Thm} \label{thm:homotopy = inner}
$vect(\alpha) \cong vect(\beta)$ if and only if $\alpha = Ad\;U\circ\beta$ 
for some unitary $U$.
\end{Thm}
We give two proofs of this theorem. The first proves the theorem directly, 
using the preceeding observation. The second proof uses a modified argument 
from \cite{Powers} (Theorem 2.4) to show that equivalence up to inner 
automorphism coincides with equivalence up to homotopy.

\textbf{Proof 1:} If $\alpha = \gamma \circ \beta$ where $\gamma$ is an 
inner automorphism, then by theorem \ref{thm:product formula} we have 
$vect(\alpha) \cong vect(\gamma) \otimes vect(\beta)$. But $vect(\gamma)$ is 
the trivial line bundle, so $vect(\gamma) \otimes vect(\beta) \cong 
vect(\beta)$.

For the converse, $vect(\alpha) \cong vect(\beta)$ means that there is an 
isomorphism between them, i.e. a unitary element in 
$Hom(vect(\alpha),vect(\beta))$, which by observation \ref{obs:hom bundle} 
can be viewed as a sub-bundle of $X \times \B(\HH)$. So we know that there 
is a unitary section $U$ of this sub-bundle (and in particular, 
$U\in\M(\cxk)$), and any such section can immediately be seen to satisfy 
$U\alpha(A)U^* = \beta(A)$ for all $A\in\cxk$, i.e. $\alpha = 
Ad\;U\circ\beta$, as required. \qed

\textbf{Proof 2:}  Suppose $U$ implements an equivalence between $\alpha$ 
and $\beta$. $U$ can be thought of as a function $u:X \arrow \UH$. We know 
that $\UH$ is contractible (by theorem \ref{thm:contractibility} for 
example, observing that $\UH \cong {\cal C}_1$). Therefore, there is a path 
of functions $u_t:X \arrow \UH$, continuous when viewed as functions $X 
\times [0,1] \arrow \UH$, such that $u_0 = u$, $u_1 \equiv 1$. Let $U_t$ be 
the corresponding path of unitaries in $\M(\cxk)$, then $\alpha_t = 
AdU_t\circ\alpha$ implements a homotopy between $\alpha$ and $\beta$.

Conversely, suppose $\alpha$ and $\beta$ are homotopic. Let $e_{ij}$ $(i,j 
\in \N)$ be matrix units for $\K$, and let $E_{ij} = 1 \otimes e_{ij} \in 
C_b(X) \otimes \K \subseteq \M(\cxk)$. So the projections $\alpha(E_{11})$, 
$\beta(E_{11})$ are homotopic (where by that we mean the canonical 
extensions of $\alpha$, $\beta$ applied to $E_{11}$). Since homotopic 
projections are unitarily equivalent (see, for example, \cite{Blackadar} 
chapter 4), we can find a unitary $V$ such that $V\alpha(E_{11})V^* = 
\beta(E_{11})$, and therefore we can find a partial isometry $W \in 
\M(\cxk)$ with $W^*W =\alpha(E_{11})$, $WW^* = \beta(E_{11})$. Let $U = 
\sum_{i=1}^{\infty}\beta(E_{i1})W\alpha(E_{1i})$, the sum taken in the 
strict topology. It is straightforward to verify that $U$ is unitary. Now, 
the span of the elements of the form $f \otimes e_{ij}$, where $f \in 
C_0(X)$ is dense in $\cxk$, and
$$
U \alpha(f \otimes e_{kl}) U^*
=
\sum_{i,j} \beta(E_{i1}) W \alpha(E_{1i}) \alpha(f \otimes e_{kl}) 
\alpha(E_{j1}) W^* \beta(E_{1j})
=
$$ $$ =
\beta(E_{k1}) W \alpha(f \otimes e_{11}) W^* \beta(E_{1l})
=
\beta(E_{k1}) W f \cdot \alpha(E_{11}) W^* \beta(E_{1l})
=
$$ $$ =
f \cdot \beta(E_{k1}) W W^* W W^* \beta(E_{1l})
=
f \cdot \beta(E_{k1}) \beta(E_{11}) \beta(E_{1l})
=
$$ $$ =
f \cdot \beta(E_{kl})
=
\beta(f \otimes e_{kl})
$$ \qed

\begin{Rmk} Notice that $\alpha$ is an automorphism if and only if 
$vect(\alpha)$ is a line bundle. So, if we restrict our attention to 
automorphisms, we get a map from the group of spectrum preserving 
automorphisms onto the group of Hermitian line bundles, whose kernel is the 
inner automorphisms. The group of line bundles is isomorphic to $H^2(X;\Z)$ 
(at least in the case in which $X$ is paracompact; see for example \cite{Karoubi}, I.3). So, we recover the exact 
sequence of Phillips and Raeburn mentioned in the introduction.
\end{Rmk}

Composition of \cxendo s maps under $vect$ to tensor products. One might ask 
if there is a natural operation on \cxendo s which is taken to the direct 
(Whitney) sum. We can indeed construct such an operation, as follows. Let 
$s_1,s_2 \in \M(\cxk)$ be two isometries such that $s_1s_1^* + s_2s_2^* = 
1$, and let $\alpha$,$\beta$ be \cxendo s of $\cxk$, then we can define 
$\alpha \oplus_{(s_1,s_2)} \beta$ by
$$
(\alpha \oplus_{(s_1,s_2)} \beta)(A) = s_1 \alpha(A) s_1^* + s_2 \beta(A) 
s_2^*
$$
While the definition depends, of course, on the choice of $s_1,s_2$, this 
direct sum construction is unique up to homotopy, and hence by theorem 
\ref{thm:homotopy = inner}, up to equivalence via an inner automorphism. The 
uniqueness up to homotopy follows easily from the fact that each pair 
$(s_1,s_2)$ corresponds in the obvious way to a map $X \arrow {\cal C}_2$, 
and any two such maps are homotopic, since ${\cal C}_2$ is contractible 
(theorem \ref{thm:contractibility}).

Define a map
$$
\Psi:vect(\alpha) \oplus vect(\beta) \arrow vect(\alpha \oplus_{(s_1,s_2)} 
\beta)
$$
by
$$
\Psi(x, T_1 \oplus T_2) = (x,s_1(x)T_1 + s_2(x)T_2)
$$
where $s_1(x)$,$s_2(x)$ are the two isometries obtain by evaluating the map 
$X \arrow {\cal C}_2$ at $x$.
To verify that this is a bundle map, we first have to see that $\Psi(x, T_1 
\oplus T_2) \in vect(\alpha \oplus_{(s_1,s_2)} \beta)$. Let us denote the 
maps $X \arrow End(\K)$ corresponding to $\alpha$,$\beta$ by $x \mapsto 
\alpha_x$,$\beta_x$, respectively, and fix $x\in X$, then for any $A\in\K$,
$$
\Psi(T_1 \oplus T_2)A
=
s_1(x)T_1A + s_2(x)T_2A
=
s_1(x)\alpha_x(A)T_1 + s_2(x) \beta_x(A)T_2
=
$$ $$
=
(\alpha \oplus_{(s_1,s_2)} \beta)_x(A)(s_1(x)T_1 + s_2(x)T_2)
=
(\alpha \oplus_{(s_1,s_2)} \beta)_x(A)\Psi(T_1 \oplus T_2)
$$
so indeed $\Psi(x, T_1 \oplus T_2) \in vect(\alpha \oplus_{(s_1,s_2)} 
\beta)$. $\Psi$ is clearly a bundle map, so now it remains to show that it 
is an isomorphism. Since $s_1(x)$,$s_2(x)$ are isometries, $\Psi$ is 
injective. For surjectivity, let $(x,T) \in vect(\alpha \oplus_{(s_1,s_2)} 
\beta)$. Let $T_1 = s_1(x)^*T$,$T_2 = s_2(x)^*T$ then for any $A\in\K$, we 
have
$$
T_1 A
=
s_1(x)^*T A
=
s_1(x)^* (\alpha \oplus_{(s_1,s_2)} \beta)_x(A) T
=
$$ $$
=
s_1(x)^* (s_1(x)\alpha_x(A)s_1^* + s_2(x) \beta_x(A)s_2^*) T
=
\alpha_x(A)s_1^*T
=
\alpha_x(A) T_1
$$
and similarly we have $T_2A = \beta_x(A) T_2$, so $(x,T_1 \oplus T_2) \in 
vect(\alpha) \oplus vect(\beta)$ and $\Psi(x,T_1 \oplus T_2) = (x,T)$, so 
$\Psi$ is surjective.

\

We conclude the section with a summary of the main results.
\paragraph{Conclusion:} To each \cxendo~ $\alpha$ of $\cxk$ we associate the 
bundle $vect(\alpha)$ of intertwining operators.
\begin{enumerate}
\item $vect(\alpha)$ is numerable, and any numerable vector 
bundle arises as $vect(\alpha)$ for some \cxendo~ $\alpha$.
\item $vect(\alpha) \cong vect(\beta)$ if and only if $\alpha = AdU \circ 
\beta$ for some $U \in \M(\cxk)$.
\item $vect(\alpha \circ \beta) = vect(\alpha) \otimes vect(\beta)$
\item $vect(\alpha)$ is a trivial bundle if and only if the endomorphism is 
inner, in the sense that there exist isometries $S_1,...,S_n \in \M(\cxk)$ 
such that $\alpha(A) = \sum S_i A S_i^*$.
\end{enumerate}

\begin{Rmk} As in the previous section, we remark that the analogous results 
for real $\Cs$-algebras can be obtained in a straightforward way, by going 
through the proofs, replacing the complex terminology by the real 
terminology.
\end{Rmk}

\section{$C_0(X)$-Endomorphisms of Stable Continuous-Trace 
$\Cs$-algebras}

Let $E$ be a numerable fiber bundle over $X$ with fiber $\K$ and 
structure group $Aut(\K)$, and let $\GE$ be the $\Cs$-algebra of sections. 
The purpose of this section is to generalize the main results of the 
previous section to this case.

Note that $\M(\GE))$ can be identified with the sections of a bundle with 
fiber $\B(\HH)$ and the same Dixmier-Douady class (see \cite{PR}, 
Proposition 2.15 and the following remark). We refer to this bundle 
henceforth as the \emph{multiplier bundle of $E$}. Denote by $\B_x$ the 
fiber of the multiplier bundle over a point $x\in X$. For a \cxendo~ 
$\alpha$, we denote by $\alpha_x$ the action of $\alpha$ on the fiber over 
$x$ (coming from the induced action of $\alpha$ on the quotient algebra 
corresponding to the primitive ideal represented by $x$), and then as 
before, we can define
$$
vect(\alpha) = \{(x,T) \; | \; T \in \B_x \; , \; TA = \alpha_x(A)T \; \; 
\forall A\in E_x\}
$$
Since $E$ is locally trivial, we see that the restriction of $vect(\alpha)$ 
to each trivialization of $E$ is locally trivial, and hence it is also 
locally trivial.

$vect(\alpha)$ determines $\alpha$ as above: if $N \subseteq X$ is such that 
$vect(\alpha)|_N$,$E|_N$ are trivial, then choose orthonormal sections 
$S_1,...,S_n$ of the restricted bundle, and then $\alpha|_N(A) = \sum 
S_iAS_i^*$ for $A \in \Gamma(E|_N)$; those restrictions determine $\alpha$.

If there is a unitary $U\in\M(\GE)$ such that $\alpha(A) = U^*\beta(A)U$ for 
all $A\in\GE$ then $TA_x = \alpha_x(A_x)T$ if and only if $(U_xT)A_x = 
\beta_x(A_x)(U_xT)$, so left multiplication by $U$ gives a vector bundle 
isomorphism $vect(\alpha) \cong vect(\beta)$.

Conversely, suppose $vect(\alpha) \cong vect(\beta)$. Note that as before, 
$Hom(vect(\alpha),vect(\beta))$ is isomorphic to the sub-bundle of the 
multiplier bundle whose fiber over $x$ is $span\{TS^* \; | \; T \in 
vect(\beta) \; , \; S \in vect(\alpha)\}$. Since there exists a unitary 
element in  $Hom(vect(\alpha),vect(\beta))$, we know that there is a unitary 
section $U$ of the above sub-bundle (and in particular, $U\in\M(\GE)$), and 
any such $U$ satisfies $U\alpha(A)U^* = \beta(A)$ for all $A\in\GE$.

So, as in the case of \cxendo s of $\cxk$, we constructed a map $vect$ from 
\cxendo s to vector bundles such that $vect(\alpha) \cong vect(\beta)$ if 
and only if $\alpha$ and $\beta$ are equivalent via an inner automorphism. 
The product formula from theorem \ref{thm:product formula} extends as well 
to this situation, with essentially the same proof.

If $E_1$,$E_2$ are two locally trivial $\K$-bundles over $X$, then $E_1 
\otimes E_2$ is a locally trivial $\K$-bundle over $X$, and if 
$\alpha_1$,$\alpha_2$ are \cxendo s of $\Gamma(E_1)$,$\Gamma(E_2)$ 
respectively then we can have a \cxendo ~ $\alpha_1 \otimes \alpha_2$ of 
$\Gamma(E_1 \otimes E_2)$.

\begin{Lemma} Let $E_1$,$E_2$ be two locally trivial $\K$-bundles over $X$, 
and let $\alpha_1$,$\alpha_2$ be \cxendo s of $\Gamma(E_1)$,$\Gamma(E_2)$, 
respectively.
\begin{enumerate}
\item $vect(\alpha\otx\beta) \cong vect(\alpha) \otimes vect(\beta)$
\item If $\varphi: \Gamma(E_1) \arrow \Gamma(E_2)$ is an isomorphism 
commuting with the multiplier action of $C_0(X)$ then 
$vect(\varphi\alpha\varphi^{-1}) \cong \nolinebreak vect(\alpha)$
\end{enumerate}
(cf \cite{RW} lemma 6.6)
\end{Lemma}
\Proof
Let $M_1$,$M_2$ denote the multiplier bundles of $E_1$,$E_2$ respectively. 
We denote by $M_1 \otimes M_2$ the fiber-wise tensor product bundle, where 
the tensor product is taken to be the spatial tensor product, and observe 
that $M_1 \otimes M_2$ is canonically the multiplier bundle of $E_1 \otimes 
E_2$.
$vect(\alpha\otimes\beta)$ is a sub-bundle of $M_1 \otimes M_2$, as is 
$vect(\alpha) \otimes vect(\beta)$, and it is immediate to verify that they 
are in fact the same. That proves $(1)$. $(2)$ follows immediately from the 
product formula (the extension of theorem \ref{thm:product formula}). \qed

We can now obtain surjectivity of $vect$ as follows (cf \cite{RW}, p. 160).
Suppose $v$ is a given numerable vector bundle. Let $E$ be our given 
$\K$-bundle, and let $E_0$ be the trivial bundle. By the results of the 
previous section, there exists a \cxendo~ $\alpha_0$ of $\Gamma(E_0) \cong 
\cxk$ with $vect(\alpha_0) \cong v$. Note that $vect(id_{\GE})$ is the 
trivial line bundle $\epsilon_0$ over $X$. Define a \cxendo~ $\alpha$ of 
$\GE$ by $\alpha = \alpha_0 \otimes id_{\GE}$, then
$$
vect(\alpha) \cong vect(\alpha_0) \otimes vect(id_{\GE}) \cong v \otimes 
\epsilon_0 \cong v
$$

If $\alpha$ is a given endomorphism $\GE$, let $E^{op}$ be the opposite bundle, and let $\alpha_0=\alpha\otimes id_{E^{op}}$, then $\alpha_0$ is an endomorphism of $E_0 \cong E \otimes E^{op}$, and as above, $vect(\alpha) \cong vect(\alpha_0)$, and we know that $vect(\alpha_0)$ is numerable, hence $vect(\alpha)$ is numerable as well.

Ilan Hirshberg

Department of Mathematics

University of California at Berkeley

Berkeley, CA 94720

USA

\

email: \emph{ilan@math.berkeley.edu}

\end{document}